\documentclass{article}

\usepackage{amsmath,amssymb,amsthm,enumerate,mathrsfs}
\usepackage[numbers]{natbib}
\usepackage{booktabs}
\usepackage{adjustbox}
\usepackage{subfigure,multirow}
\usepackage{tikz}
\usepackage{appendix}
\usepackage{authblk}
\usepackage{geometry}
\usepackage[ruled,linesnumbered]{algorithm2e}
\usepackage{algpseudocode}%
\usepackage{longtable}
\usepackage{caption}
\usepackage{lipsum}
\usepackage{hyperref}
\usepackage{footnote}
\usepackage[misc]{ifsym}

\newcommand{\T}{\text{T}}
\newcommand{\M}{Markowitz }
\newcommand{\BB}{Barzilai-Borwein}

\renewcommand{\thefootnote}{\fnsymbol{footnote}}

\theoremstyle{plain}

\theoremstyle{remark}
\newtheorem{remark}{Remark}

\newtheorem{property}{Property}

\newcommand\blfootnote[1]{% 
	\begingroup 
	\renewcommand\thefootnote{}\footnote{#1} 
	\addtocounter{footnote}{-1}% 
	\endgroup 
}

\title{An adaptive ADMM with regularized spectral penalty for sparse portfolio selection}
\author{Xin Xu\textsuperscript{\rm *}}
\date{\today}

\begin{document}
\maketitle

\blfootnote{\Letter \text{Xin Xu} \par \; \text{xustonexin@gmail.com}} 
\footnotetext[1]{School of Mathematics, Southwestern University of Finance and Economics, Chengdu 611130, Sichuan, China}

\begin{abstract}
	The mean-variance (MV) model is the core of modern portfolio theory. Nevertheless, it suffers from the over-fitting problem due to the estimation errors of model parameters. We consider the $\ell_{1}$ regularized MV model, which adds an $\ell_{1}$ regularization term in the objective to prevent over-fitting and promote sparsity of solutions. By investigating the relationship between sample size and over-fitting, we propose an initial regularization parameter scheme in the $\ell_{1}$ regularized MV model. Then we propose an adaptive parameter tuning strategy to control the amount of short sales. ADMM is a well established algorithm whose performance is affected by the penalty parameter. In this paper, a penalty parameter scheme based on regularized \BB \ step size is proposed, and the modified ADMM is used to solve the $\ell_{1}$ regularized MV problem.  Numerical results verify the effectiveness of the two types of parameters proposed in this paper.  
	
	
\end{abstract}

\section{Introduction}\label{sec:introduction}
Portfolio theory address the decision problem of how to allocate finite resources among several competing assets. The Mean-Variance (MV) model was proposed by \M \cite{Markowitz1952PortfolioSelection} and is the core of the modern portfolio theory. This model provides a trade-off between expected return and risk, allowing investors to obtain a given expected return with minimal risk. In subsequent work \cite{Markowitz1959PortfolioSelection}, \M enhanced his theory and argued that under mild conditions, the portfolio derived from the mean-variance efficient frontier will approach maximizing the investor's expected utility. However, it has been shown that the MV is usually suboptimal with poor out-of-sample performance due to its high sensitivity to estimation errors of covariance matrices and expected returns \cite{Best1991SensitivityMeanVariance,Chopra1993EffectErrorsMeans,BrittenJones1999SamplingErrorEstimates}. Generally, it is more difficult to estimate means than covariances of assets returns and also that the estimation errors of means has more influence on portfolio than that of covariances \cite{Merton1980estimatingexpectedreturn}.

To deal with this problem, \cite{DeMiguel2009GeneralizedApproachPortfolio} provided a general framework that relies on solving the traditional minimum-variance model but subject to the additional constraint that the norm of the portfolio-weight vector be smaller than a given threshold. The regularization method adopts the idea of Lagrangian function and adds this constraint to the objective. It can avoid over-fitting and thus improve out-of-sample performance \cite{Fastrich2014Constructingoptimalsparse,Gunjan2022briefreviewportfolio}, and is an effective method to deal with the estimation error in portfolio selection. Originating from the distinct features of $\ell_{0}$, $\ell_{1}$, and $\ell_{2}$ norms, constraints based on these norms frequently arise in portfolio optimization. It is worth noting that the $\ell_{0}$ defines a semi-norm. 

There are some efforts on these constrained MV models. In \cite{Cesarone2012newmethodmean,DiLorenzo2012concaveoptimizationbased,Wang2023l0normbased}, $\ell_{0}$ norm constrained MV have been proposed to control the cardinality of the  portfolio, which are solved using concave minimization or regularization techniques. Nevertheless, $\ell_{0}$ norm optimization problems are NP-hard. To circumvent this issue, an alternative solution is relaxing $\ell_{0}$ to the $\ell_{1}$ norm, corresponding to the LASSO approach, which results in a tractable convex optimization problem. It is worth noting that, with relaxation, the target cardinality behaves as a sparsity control parameter. \cite{Brodie2009SparsestableMarkowitza} used $\ell_{1}$ regularization to obtain a sparse and stable optimal portfolio, and proved that imposing a positive constraint on the portfolio weights is equivalent to penalizing a largest $\ell_{1}$ regularization to the objective function. \cite{Fan2012VastPortfolioSelection} provided theoretical and empirical evidence that using $\ell_{1}$ regularization to limit the total number of investments can achieve a sparse and stable portfolio without causing the accumulation of estimation errors. The common consensus is that the $\ell_{2}$ norm constraint stabilizes the computation by improving the condition number of the problem resulting in strong out-of-sample performance. \cite{Yen2014Solvingnormconstrained} considered a weighted $\ell_{1}$ and $\ell_{2}$ square norm penalty. \cite{Fastrich2014Constructingoptimalsparse} proposed a class of regularization operators including convex and non-convex penalties. \cite{Corsaro2018Adaptive$$l_1$$l} proposed an $\ell_{1}$ regularized MV model with an adaptive regularization parameter and used the Bregman iteration technique to solve it. \cite{Zhao2021Optimalportfolioselections,Wu2024Sparseportfoliooptimization} combined the advantages of $\ell_{1}$ and $\ell_{2}$ regularization to improve the sparsity and stability of the solution. Since portfolio determines the amount of capital to be invested in each available security, sparsity means that the funds will be invested in a small number of securities, the active positions. This allows investors to reduce the number of positions to monitor and transaction costs, especially for small investors.

The alternating direction method of multipliers (ADMM) is a well established method for solving the minimization problem with linear constrains and a separable objective function. It has been successfully applied in different fields, for example, the robust principal component analysis model with noisy and incomplete data \cite{Tao2011RecoveringLowRank}, the nonnegative matrix complete \cite{Xu2012alternatingdirectionalgorithm}, the image restoration \cite{Zietlow2024ADMMTGVimage}. The ADMM decomposes the original problem into several subproblems. The resultant subproblems can be solved efficiently, especially when they have closed-form expressions. It is worth noting that \cite{Shi2024CardinalityConstrainedPortfolioa} proposed a cardinality constrained portfolio optimization via ADMM.

This paper focuses on an adaptive ADMM algorithm for an $\ell_{1}$ regularized mean-variance portfolio model. Specifically, we first propose a regularized \BB \ step size (spectral) penalty parameter for ADMM, which integrates the residual balancing \cite{He2000AlternatingDirectionMethod} and spectral penalty \cite{Xu2017AdaptiveADMMSpectral} strategies, and fully utilizes the convergence characteristics of ADMM. As a means of preventing over-fitting, regularization balances the estimation and generalization capabilities of the model. When estimating a model on a small dataset, one often suffers from the over-fitting problem, and a high regularization is desired. In another extreme case, if the training dataset is large enough, the estimated model is close to the truth, and a small regularization is reasonable. In addition, in the $\ell_{1}$ regularized model with the constraint $x^{\T}\mathbf{1}_{n}=1$, as the number of assets increases, the solutions present a sparse distribution naturally. Based on these observations, we propose an initial regularization parameter scheme that employs the sample size and the number of assets. With the control of short sales as the financial goal, we adjust the regularization parameter when the number of short sales exceeds a preset value. 

This paper is organized as follows. In Section \ref{sec:Background}, we present the $\ell_{1}$ regularized portfolio model and review the basic framework of ADMM as well as its two effective penalty parameter methods. In Section \ref{sec:RBB}, we first propose a regularized spectral penalty scheme for ADMM, which combines the previous two penalty parameters, and then provide a regularization parameter strategy of the $\ell_{1}$ regularized mean-variance portfolio model, which is related to the sample size and the  number of assets. In Section \ref{numerical}, numerical experiments are presented. 
 
\section{Background}\label{sec:Background}
In this section, we review some basic content, including the mean-variance portfolio model, the ADMM framework, and the selection of penalty parameters in ADMM.  
\subsection{Portfolio selection model}\label{sec:Portfolio}
We refer to the classical \M mean-variance framework. Given $n$ traded assets, the core of the problem is to establish the amount of capital to be invested in each available security.
 
Specially, it is assume that one unit of capital is available. A portfolio is defined as a weight vector $x=(x_{1},x_{2},\ldots,x_{n})^{\T}$, where $x_{i}$ represents the proportion of capital to be invested in the $i$-th security. Asset return are assume to be stationary. If we denote with $\mu=(\mu_{1},\mu_{1},\ldots,\mu_{n})^{\T}$ the expected asset return, then the expected portfolio return is 
\begin{equation}\label{expect}
	\sum_{i=1}^{n}x_{i}\mu_{i}.
\end{equation}  
We denote with $\sigma_{ij}$ the covariance between returns of securities $i$ and $j$. The portfolio risk is measured by 
\begin{equation*}
	V=\sum_{i=1}^{n}\sum_{j=1}^{n}\sigma_{ij}x_{i}x_{j}.
\end{equation*}
Let $e$ be the fixed expected portfolio return and $C$ the covariance matrix of returns. We assume $C$ to be positive definite. Portfolio selection is formulated as the following quadratic programming problem:
\begin{equation}\label{mv}
	\begin{aligned}
	\mathop{\text{min}}\limits_{x\in\mathbb{R}^{n}}\quad \frac{1}{2}x^{\T}Cx\quad
	\text{subject to} \quad x^{\T}\mu=e, \
	     x^{\T}\mathbf{1}_{n}=1,
	\end{aligned}
\end{equation} 
where $\mathbf{1}_{n}$ is the vector of ones of length $n$. The first constraint fixes the expected return. The second one is a budget constraint which establishes that all the available capital is invested. In this paper, we consider the following $\ell_{1}$ regularized problem:
\begin{equation}\label{reguMV}
\begin{aligned}
\mathop{\text{min}}\limits_{x\in\mathbb{R}^{n}}\quad \frac{1}{2}x^{\T}Cx+\lambda\|x\|_{1}\quad
\text{subject to} \quad x^{\T}\mu=e, \ x^{\T}\mathbf{1}_{n}=1,
\end{aligned}
\end{equation} 
where $\lambda$ is the regularization parameter. Problem \eqref{reguMV} can be formulated as the constrained nonlinear optimization problem as follows
\begin{equation}\label{portfolio}
\begin{aligned}
\mathop{\text{min}}\limits_{x\in\mathbb{R}^{n}}\quad \frac{1}{2}x^{\T}Cx+\lambda\|x\|_{1}, \ 
\text{subject to} \quad Dx=b, 
\end{aligned}
\end{equation}
where 
$
D=\begin{pmatrix}
\mu^{\T}\\
\mathbf{1}_{n}^{\T}
\end{pmatrix}
$
and
$
b=\begin{pmatrix}
e\\
1
\end{pmatrix}
.
$
Now we reformulate Problem \eqref{portfolio} as
\begin{equation}\label{ADMMportfolio}
\begin{aligned}
\mathop{\text{min}}\limits_{x\in\mathbb{R}^{n},z\in\mathbb{R}^{n}}\quad f(x)+g(z),\quad 
\text{subject to} \quad x-z=0,
\end{aligned}
\end{equation}
where $f(x)=\frac{1}{2}x^{\T}Cx$ restricted to the set $\{x|Dx=b\}$ and $g(z)=\lambda\|z\|_{1}$.

\subsection{ADMM framework}\label{ADMMRBB}

ADMM is an important algorithm in the optimization community, which solves optimization problems of the from 
\begin{equation}\label{original}
\begin{aligned}
\mathop{\text{min}}\limits_{x\in\mathbb{R}^{n},z\in\mathbb{R}^{m}} \quad & f(x)+g(z)\\
\text{subject to} \quad & Ax+Bz=c,
\end{aligned}
\end{equation}
where $f:\mathbb{R}^{n}\rightarrow\mathbb{R}$ and $g:\mathbb{R}^{m}\rightarrow\mathbb{R}$ are closed, proper, and convex functions, matrices $A\in\mathbb{R}^{p\times n}$, $B\in\mathbb{R}^{p\times m}$, vector $c\in\mathbb{R}^{p}$. Note that the objective function in \eqref{ADMMportfolio} is separable and the constraint is liner, which is exactly the form that ADMM can solve with $A=I$, $B=-I$ and $c=0$. The augmented Lagrangian of problem \eqref{original} is
\begin{equation}\label{Lagrangian}
	L_{\rho}(x,z,y)=f(x)+g(z)+y^{\T}(c-Ax-Bz)+\frac{\rho}{2}\|c-Ax-Bz\|_{2}^{2},
\end{equation}
where $y\in\mathbb{R}^{p}$ is the dual variable (Lagrange multiplier), $\rho>0$ is the penalty parameter. Instead of jointly solving for the primal variables $x$ and $z$, ADMM updates $x$ and $z$ alternately. Starting from an initialization ($x^{0}, z^{0}, y^{0}$) and $\rho>0$, at each iteration, in a scaled form, ADMM updates the primal and dual variables as 
	\begin{align}
	x^{k+1}=&\mathop{\text{argmin}}\limits_{x}f(x)+\frac{\rho}{2}\|c-Ax-Bz^k+\frac{y^k}{\rho}\|_{2}^{2}\label{ADMM1}\\
	z^{k+1}=&\mathop{\text{argmin}}\limits_{z}g(z)+\frac{\rho}{2}\|c-Ax^{k+1}-Bz+\frac{y^k}{\rho}\|_{2}^{2}\label{ADMM2}\\
	y^{k+1}=&y^k+\rho(c-Ax^{k+1}-Bz^{k+1})\label{ADMM3},
	\end{align}
where $k$ is iteration number. In the iterative format of ADMM, $\rho$ is the only freely optional one. The convergence of ADMM can be monitored using primal and dual residuals, both of which approach zero as the iterates approximate optimal solution, and which are defined as 
\begin{equation}\label{residuals}
	r^k=c-Ax^k-Bz^k, \quad \text{and}\quad d^k=\rho A^{\T}B(z^k-z^{k-1}),
\end{equation}
respectively \cite[p.18]{S.Boyd2011DistributedOptimizationStatistical}. The iteration is stopped when $\|r^k\|_{2}\le\varepsilon^{tol}\max\{\|Ax^k\|_2, \ \|Bz^k\|_{2}, \ \|c\|_2\}$ and $\|d^k\|_2\le\varepsilon^{tol}\|A^{\T}y^{k}\|_{2}$, where $\varepsilon^{tol}$ is the stopping tolerance. For any fixed value of positive $\rho$, ADMM can be shown to converge to a solution of problem \eqref{original} under mild assumptions (see, eg.,  \cite{Gabay1983ApplicationsMethodMultipliers,Eckstein1992DouglasRachfordsplittinga,S.Boyd2011DistributedOptimizationStatistical,Deng2015GlobalLinearConvergence}). \cite{He2012$O1/n$ConvergenceRate} proved that the convergence rate of ADMM is $O(1/k)$ in the ergodic sense in 2012 via variational inequality for convex problems. However, it is well know that in practice, the rate of convergence of ADMM is closely related to the size of penalty parameter $\rho$. For non-expert users, it is difficult to get a suitable value of $\rho$.

\subsubsection{Review two penalty parameter methods}
The residual balancing strategy and the spectral gradient step size scheme are two effective penalty parameter methods in practice.

The residual balancing (RB) \cite{He2000AlternatingDirectionMethod} strategy is derived from the observation on the convergence behavior of ADMM in  \cite{He2000AlternatingDirectionMethod}, where the distance to convergence is defined as  $D^{k+1}=\|r^{k+1}\|_{2}^{2}+\|d^{k+1}\|_{2}^{2}$. As its name suggested, RB balances the convergence progress of $\|r^{k}\|_{2}$ and $\|d^{k}\|_{2}$, keeping them at similar magnitude, and achieves simultaneous convergence at the optimal solution. Specifically, the form of RB is as 
\begin{equation}\label{RB}
	\rho^{k+1}=\begin{cases}
	 \eta\rho^{k} \quad & \text{if}\quad\|d^{k}\|_{2}>\mu\|r^{k}\|_{2}\\
	 \rho^{k}/\eta \quad & \text{if}\quad\|r^{k}\|_{2}>\mu\|d^{k}\|_{2}\\
	 \rho^{k} \quad & \text{otherwise},
	\end{cases}
\end{equation}
for parameters $\eta>1$ and $\mu>1$, which is based on the observation that increasing $\rho^{k}$ strengthens the penalty term, resulting in smaller primal residuals but larger dual ones; conversely, decreasing $\rho^{k}$ leading to smaller dual residuals but larger primal ones. But one drawback of RB is sensitive to problem scaling \cite{Mccann2024RobustSimpleADMM}. In addition, the performance of RB is affected by the initial $\rho^{0}$ and the values of $\eta, \mu$. If these values are not set appropriately, ADMM may not converge before $\rho^{k}$ is fixed. 

Adaptive spectral gradient step size penalty (SP) strategy was introduced in \cite{Xu2017AdaptiveADMMSpectral}, which exploits the close relationship between ADMM and Douglas-Rachford (DR) splitting applied to the dual problem of \eqref{original} and transforms the selection of the penalty parameter $\rho^{k}$ in ADMM to the selection of the step size in DR. We now give a brief review.

Given a closed convex (proper) function $h$ defined on $\mathbb{R}^{n}$, the Fenchel conjugate of $h$ is the closed convex (proper) function $h^{*}$ defined by
$$h^{*}(x^{*})=\mathop{\text{sup}}\limits_{x}\{x^{\T}x^{*}-h(x)\}=-\mathop{\text{inf}}\limits_{x}\{h(x)-x^{\T}x^{*}\},$$
(see \cite[p.104]{Rockafellar.1970ConvexAnalysis}).
The dual problem of \eqref{original} (cf. \cite[\S2.3]{Esser2009ApplicationsLagrangianBased}) is given by 
\begin{equation}\label{dualmin}
\mathop{\text{min}}\limits_{y\in\mathbb{R}^{p}} f^{*}(A^{\T}y)-y^{\T}c+g^{*}(B^{\T}y).
\end{equation}
where $f^*$ and $g^*$ denote the Fenchel conjugate of $f$ and $g$, respectively. 
Note that this is an unconstrained optimization problem with respect to $y$. Define operators $\Psi$ and $\Phi$ by
\begin{align}
\Psi(y)&=A\partial f^{*}(A^{\T}y)-c,\\
\Phi(y)&=B\partial g^{*}(B^{\T}y).
\end{align}
It can be proved that solving \eqref{original} by ADMM is equivalent to solving the dual problem \eqref{dualmin} by means of the DR method \cite[\S3.2.3]{Esser2009ApplicationsLagrangianBased}, which is equivalent to applying the DR scheme to 
\begin{equation}\label{DRSS}
0\in \Psi(y) + \Phi(y).	
\end{equation}
By formally applying DR splitting to \eqref{DRSS} with $\rho^{k}$, two sequences $\{\bar{y}^{k}\}$ and $\{y^{k}\}$ are generated such that
\begin{align}
0&\in\frac{\bar{y}^{k+1}-y^{k}}{\rho^{k}}+\Psi(\bar{y}^{k+1})+\Phi(y^{k}),\label{DRS1}\\
0&\in\frac{y^{k+1}-y^{k}}{\rho^{k}}+\Psi(\bar{y}^{k+1})+\Phi(y^{k+1}).\label{DRS2}
\end{align}
Then, Proposition $1$ in \cite{Xu2017AdaptiveADMMSpectral} proves that the choice of the parameter $\rho^{k}$ that guarantees the minimal residual of $\Psi(y^{k+1})+\Phi(y^{k+1})$ in DR steps is given by 
\begin{equation}\label{ab}
\rho^{k}=1/\sqrt{\alpha_{k}\beta_{k}},
\end{equation}
where $1/\alpha_{k}, 1/\beta_{k}>0$ are \BB\cite{Barzilai1988TwoPointStep}(BB) step sizes deriving from the following least squares model
\begin{align}
\alpha_{k}&=\mathop{\text{argmin}}\limits_{\alpha\in\mathbb{R}} \|\alpha(\bar{y}^{k}-\bar{y}^{k-1})-(\psi^{k}-\psi^{k-1})\|_{2},\label{B1}\\
\beta_{k}&=\mathop{\text{argmin}}\limits_{\beta\in\mathbb{R}} \|\beta(y^{k}-y^{k-1})-(\phi^{k}-\phi^{k-1})\|_{2},\label{B11}
\end{align}
where $\psi^{k}\in\Psi(\bar{y}^{k})$, $\psi^{k-1}\in\Psi(\bar{y}^{k-1})$, $\phi^{k}\in\Phi(y^{k})$, and  $\phi^{k-1}\in\Phi(y^{k-1})$. 

In practice, we only compute $\rho^{k}$ using the ADMM steps without considering the dual problem \eqref{dualmin}, due to the theoretical connection between the primal and dual variables. Specifically, investigating the relationship between the optimality conditions of the ADMM steps \eqref{ADMM1},\eqref{ADMM2} and the DR steps \eqref{DRS1}, \eqref{DRS2}, we obtain the form of $\bar{y}^{k}$ as 
\begin{equation}\label{bary}
\bar{y}^{k}=y^{k-1}+\rho^{k-1}(c-Ax^{k}-Bz^{k-1}),
\end{equation}
(see \cite[p.3]{Crisci2023AdaptivePenaltyParameter}).  In addition, 
\begin{equation}\label{AB}
Ax^{k}-c\in A\partial f^{*}(A^{\T}\bar{y}^{k})-c=\Psi(\bar{y}^{k}),\quad\text{and}\quad Bz^{k}\in B\partial g^{*}(B^{\T}y^{k})=\Phi(y^{k})
\end{equation}
are deduced from the optimality conditions of the ADMM steps \eqref{ADMM1}, \eqref{ADMM2},  \eqref{bary} and \cite[Corollary 23.51]{Rockafellar.1970ConvexAnalysis}, respectively. Finally, considering \eqref{bary} and \eqref{AB}, one can define $\Delta y^{k-1}=y^{k}-y^{k-1}$,  $\Delta\bar{y}^{k-1}=\bar{y}^{k}-\bar{y}^{k-1}$,    $\Delta\Psi^{k-1}=\Psi(\bar{y}^{k})-\Psi(\bar{y}^{k-1})=A(x^{k}-x^{k-1})$ and $\Delta\Phi^{k-1}=\Phi(y^{k})-\Phi(y^{k-1})=B(z^{k}-z^{k-1})$. Then the two BB1 spectral gradient scalars in \eqref{B1} and \eqref{B11} can be rewritten as
\begin{equation}\label{BB1}
\begin{aligned}
\alpha_{k}^{BB1}&=\mathop{\text{argmin}}\limits_{\alpha\in\mathbb{R}} \|\alpha\Delta\bar{y}^{k-1}-\Delta\Psi^{k-1}\|_{2}=\frac{(\Delta\bar{y}^{k-1})^{\T}\Delta\Psi^{k-1}}{\|\Delta\bar{y}^{k-1}\|_2^2},\\
\beta_{k}^{BB1}&=\mathop{\text{argmin}}\limits_{\beta\in\mathbb{R}} \|\beta\Delta y^{k-1}-\Delta\Phi^{k-1}\|_{2}=\frac{(\Delta y^{k-1})^{\T}\Delta\Phi^{k-1}}{\|\Delta y^{k-1}\|_2^2}.
\end{aligned}	
\end{equation}
Similarly, the corresponding BB2 spectral gradient scalars are 
\begin{equation}\label{BB2}
\begin{aligned}
\alpha_{k}^{BB2}&=\mathop{\text{argmin}}\limits_{\alpha\in\mathbb{R}} \|\Delta\bar{y}^{k-1}-\alpha^{-1}\Delta\Psi^{k-1}\|_{2}=\frac{(\Delta\bar{y}^{k-1})^{\T}\Delta\Psi^{k-1}}{\|\Delta\Psi^{k-1}\|_2^2},\\
\beta_{k}^{BB2}&=\mathop{\text{argmin}}\limits_{\beta\in\mathbb{R}} \|\Delta y^{k-1}-\beta^{-1}\Delta\Phi^{k-1}\|_{2}=\frac{(\Delta y^{k-1})^{\T}\Delta\Phi^{k-1}}{\|\Delta \Phi^{k-1}\|_2^2}.
\end{aligned}	
\end{equation}
The preceding quasi-Newton conditions \eqref{B1} assume the change in dual gradient is linearly proportional to the change in the dual variables. To test the validity od this assumption, we measure the correlation between these quantities \cite{Xu2017AdaptiveADMMSpectral}:
\begin{equation}\label{corr}
\alpha_{k}^{corr}=\frac{(\Delta\bar{y}^{k-1})^{\T}\Delta\Psi^{k-1}}{\|\Delta\bar{y}^{k-1}\|_{2}\|\Delta\Psi^{k-1}\|_2} \quad\text{and}\quad\beta_{k}^{corr}=\frac{(\Delta y^{k-1})^{\T}\Delta\Phi^{k-1}}{\|\Delta y^{k-1}\|_{2}\|\Delta\Phi^{k-1}\|_{2}}.
\end{equation}
If $\alpha_{k}^{corr}$ and $\beta_{k}^{corr}$ are greater than a positive threshold $\bar{\epsilon}$, the obtained step size is considered reliable. Based on this criterion, the safeguarded spectral adaptive penalty parameter rule is 
\begin{equation}\label{rho}
\rho^{k+1}=\begin{cases}
\sqrt{\hat{\alpha}_{k}\hat{\beta}_{k}}\quad&\text{if}\quad \alpha_{k}^{corr}>\bar{\epsilon}\quad\text{and}\quad\beta_{k}^{corr}>\bar{\epsilon}\\
\hat{\alpha}_{k}\quad&\text{if}\quad \alpha_{k}^{corr}>\bar{\epsilon}\quad\text{and}\quad\beta_{k}^{corr}\le\bar{\epsilon}\\
\hat{\beta}_{k}\quad&\text{if}\quad \alpha_{k}^{corr}\le\bar{\epsilon}\quad\text{and}\quad\beta_{k}^{corr}>\bar{\epsilon}\\
\rho^{k}\quad&\text{otherwise},
\end{cases}
\end{equation} 
where $\bar{\epsilon}\in(0,1)$, $\hat{\alpha}_{k}=1/\alpha_{k}$ and $\hat{\beta}_{k}=1/\beta_{k}$ are the spectral gradient step sizes given by \eqref{BB1} or \eqref{BB2}.

We note that the rule for calculating the spectral gradient step sizes affect the performance of ADMM. \cite{Xu2017AdaptiveADMMSpectral} adopts a generalized strategy of alternating step sizes \cite{B.Zhou2006Gradientmethodsadaptive}, where a short step size is selected when $\alpha_{k}^{BB2}/\alpha_{k}^{BB1}<0.5$, otherwise a long step size is selected. In the next section, we adopt a regularized \BB \ (RBB) step size based on the properties of ADMM, and provide a more appropriate adaptive spectral gradient step size penalty parameter for ADMM.   

\section{Selecting the parameters of ADMM and regularized MV model}\label{sec:RBB}
In this section, we propose a regularized spectral \cite{an2024regularizedbarzilaiborweinmethod} penalty parameter $\rho_{k}^{RBB}$ for ADMM, which can be regard as a fusion of the two penalty parameter methods in the previous subsection. Then, we provide a regularization parameter scheme $\lambda_{k}$ of the $\ell_{1}$ regularized MV portfolio model \eqref{reguMV}.

\subsection{Regularized spectral penalty for ADMM}
Assume that there exists a symmetric positive matrix $H_{k}$ such that $\Delta\Psi^{k-1}=H_{k}\Delta y^{k-1}$. In \cite{an2024regularizedbarzilaiborweinmethod}, the authors introduce a regularized BB scalar as 
\begin{align}\label{equ:alphaRBB}
\alpha_{k}^{RBB}=\frac{(\Delta y^{k-1})^{\T}\Delta\Psi^{k-1}+\tau_{k}\|\Delta\Psi^{k-1}\|_{2}^{2}}{\|\Delta y^{k-1}\|_{2}^{2}+\tau_{k}(\Delta y^{k-1})^{\T}\Delta\Psi^{k-1}},
\end{align}
which is derived from the following regularized least squares model
\begin{equation*}
\alpha_{k}^{RBB}=\mathop{\text{argmin}}_{\alpha\in\mathbb{R}} \Big\{\Vert\alpha\Delta y^{k-1}-\Delta\Psi^{k-1}\Vert_{2}^{2} + \tau_{k}\Vert\alpha \sqrt{H_{k}}\Delta y^{k-1}- \sqrt{H_{k}}\Delta\Psi^{k-1}\Vert_{2}^{2}\Big\},
\end{equation*}
where regularized parameter $\tau_{k}\in[0,\infty)$. 

\begin{property}\cite[Them.2.2]{an2024regularizedbarzilaiborweinmethod}
	Assume that $(\Delta y^{k-1})^{\T}\Delta\Psi^{k-1}>0$, $\tau_{k}\ge0$. Then  $\alpha_{k}^{RBB}\in[\alpha_{k}^{BB1}, \  \alpha_{k}^{BB2}]$ and it is monotonically increasing with respect to the parameter $\tau_{k}$.
\end{property}
Using a similar process, we obtain $\beta_{k}^{RBB}$ as
\begin{align}\label{equ:betaRBB}
\beta_{k}^{RBB}=\frac{(\Delta y^{k-1})^{\T}\Delta\Phi^{k-1}+\tau_{k}\|\Delta\Phi^{k-1}\|_{2}^{2}}{\|\Delta y^{k-1}\|_{2}^{2}+\tau_{k}(\Delta y^{k-1})^{\T}\Delta\Phi^{k-1}}.
\end{align}
Combining the rule of penalty parameter \eqref{ab} with \eqref{equ:alphaRBB} and \eqref{equ:betaRBB}, we get the RBB step size penalty parameters in ADMM as 
\begin{equation}\label{rhoRBB}
	\rho_{k}^{RBB}=1/\sqrt{\alpha_{k}^{RBB}\beta_{k}^{RBB}}.
\end{equation}

In the spectral gradient step size methods, alternating strategy between long and short step sizes can improve the performance of algorithm. By adaptively adjusting the  parameter $\tau_{k}$, RBB method can be viewed as a continuous alternating step size method. One of its advantages is that if the regularization parameter is appropriately selected according to the nature of problem, then RBB generates reasonable step sizes, avoiding manual setting of alternating thresholds. Based on the analysis in the preceding section, the penalty parameter selection problem in ADMM is now transformed to the regularization parameter $\tau_{k}$ selection problem in RBB. 

We now consider how to select appropriate  parameter $\tau_{k}$ in the RBB step sizes \eqref{equ:alphaRBB} and \eqref{equ:betaRBB}.
The RB strategy \eqref{RB} is an effective method for adjusting the penalty parameter, which describes the convergence properties of ADMM \cite{He2000AlternatingDirectionMethod}. Combined with the spectral penalty parameter scheme \eqref{rhoRBB}, the RB strategy inspires us that if $\|r^k\|_{2}>\|d^k\|_{2}$, then in the next iteration, a small spectral gradient step size is reasonable (corresponding to a large regularization parameter $\tau_{k+1}$); otherwise, a large spectral gradient step size (corresponding to a small regularization parameter $\tau_{k+1}$) is selected when $\|r^k\|_{2}<\|d^k\|_{2}$. Based on these analyses, a natural conclusion is that the size of regularization parameter should be proportional to the ratio of $\|r^k\|_{2}$ and $\|d^k\|_{2}$. Therefore, we obtain the regularization parameter as follows
\begin{equation}\label{tauk}
	\tau_{k+1}=\Big(\frac{\|r^{k}\|_{2}}{\|d^{k}\|_{2}}\Big)^{q},
\end{equation} 
where $q>0$ is a constant that acts as a scaling factor. 
\begin{remark}
	Similar to the procedure in \cite{Mccann2024RobustSimpleADMM}, the penalty parameter generated by the RBB step size is scaling covariant and translation invariant.
\end{remark}

\subsection{Selecting regularization parameter of regularized MV model}
In this subsection, we first give the initial regularization parameter $\lambda_{0}$ based on the sample size and the number of assets, and then consider adaptively adjusting the regularization parameters based on the financial goal: controlling the number of short sales.

For analysis, the $\ell_{1}$ regularized model \eqref{reguMV} can be rewritten as 
\begin{equation}\label{sparsestable}
\begin{aligned}
\mathop{\text{min}}\limits_{x\in\mathbb{R}^{n}}\quad \frac{1}{2m}\|e\mathbf{1}_{m}-Rx\|_{2}^{2}+\lambda\|x\|_{1}, \ 
\text{subject to} \quad x^{\T}\mu=e, \ x^{\T}\mathbf{1}_{n}=1,
\end{aligned}
\end{equation}
where matrix $\mathbb{R}^{m\times n}$ is the  historical returns of asset $i$ on its $i$-th column over $m$ observation periods, see \cite{Brodie2009SparsestableMarkowitza}. The observation data contained in this matrix is used to estimate two important parameters in MV model: expected return and covariance matrix. If the number of the chosen observations $m$ is small compared to the number of assets, then the sample covariance matrix becomes ill-conditioned, one suffers from the over-fitting problem \cite{Dai2018generalizedapproachsparse,Ballal2021AdaptiveRegularizationApproach}. In this case, a high regularization promotes the sparsity of the solutions and mitigates over-fitting. On the other hand, if the number of samples is large enough, then the model obtained is truth and a small amount of regularization is reasonable. Based on these, we regard $\frac{1}{m}$ as a factor of the initial regularization parameter $\lambda_{0}$.

The second factor that affects regularization is the number of assets $n$. For any $x\neq 0$, the following inequalities hold
\begin{equation*}
	\frac{1}{\sqrt{n}}\le\frac{\|x\|_{2}}{\|x\|_{1}}\le 1.
\end{equation*}
If the constraint $x^{\T}\mathbf{1}_{n}=1$ is imposed, then the above inequalities describe the sparsity of data $x$ in high-dimensional space. In other words, as the assets size $n$ increases, the solution of the linearly constrained $\ell_{1}$ regularized model will naturally present a sparse distribution. At this point, a slight regularization has a positive effect. In short, the second factor that affects the initial regularization parameter is $\frac{1}{n}$. Combining these two factors, we obtain an initial regularization parameter as
\begin{equation}\label{initiallambda}
	\lambda_{0}=\frac{1}{mn}.
\end{equation}

From the second constraint in \eqref{sparsestable} it follows that the objective function can be written as
$$\frac{1}{2m}\|e\mathbf{1}_{m}-Rx\|_{2}^{2}+2\lambda\mathop{\sum}_{i:x_{i}<0}|x_{i}|+\lambda.$$
This implies that the $\ell_{1}$ penalty is equivalent to the penalty for short sales. Moreover, suppose that the two portfolios $x_{\lambda_{1}}$ and $x_{\lambda_{2}}$ are minimizers for the objective function in \eqref{sparsestable}, corresponding to the values $\lambda_{1}$ and $\lambda_{2}$, respectively, and both satisfy the constraint $x^{\T}\mathbf{1}_{n}=1$. We can obtain 
\begin{equation}\label{monotone}
	(\lambda_{1}-\lambda_{2})(\|x_{\lambda_{2}}\|_{1}-\|x_{\lambda_{1}}\|_{1})\ge 0,
\end{equation}
see \cite{Brodie2009SparsestableMarkowitza}. If $\lambda_{1}\ge\lambda_{2}$, it follows from \eqref{monotone} that $\|x_{\lambda_{2}}\|_{1}\ge\|x_{\lambda_{1}}\|_{1}$. If all the elements of $x_{\lambda_{1}}$ are nonnegative, but some of the $x_{i}^{\lambda_{2}}$ are negative, then we have $\|x_{\lambda_{2}}\|_{1}\ge\|x_{\lambda_{1}}\|_{1}=1$. If the regularization parameter is increased to a threshold, an all-positive portfolio is obtained. Based on this fact, if the short-position exceed $sn$, which is the threshold set by the investor, then we can adaptively adjust the regularization parameter as follows
\begin{equation}\label{lambda}
	\lambda_{k+1} = \frac{sm}{sn}\lambda_{k},
\end{equation} 
where $sm$ is the number of negative elements in $x^{k}$. In practice, only a few number of parameter adjustments are required to control the short sales, due to the existence of the first constraint in \eqref{portfolio}.

A general ADMM method with RBB step size penalty for $\ell_{1}$ regularized portfolio selection is described in Algorithm \ref{alg:AADMM}. The $\bar{n}$ in the Algorithm controls the frequency of penalty parameter update. If the value is infinite, penalty is a constant; if the value is $1$, penalty changes in each iteration. The convergence analysis of the ADMM algorithm with variable penalty was given by \cite[Them 4.1]{He2000AlternatingDirectionMethod}, under suitable conditions. Although this convergence analysis cannot straightforwardly be applied to our case, this issue may be bypassed in practice by controlling the update frequency, as observed in \cite{Xu2017AdaptiveADMMSpectral,Crisci2023AdaptivePenaltyParameter}. The proposed RBB step size penalty parameter implicitly adopts the residual balance idea and further improves the convergence of ADMM in practice. 

Specifically, in the portfolio model \eqref{ADMMportfolio}, the solution $x^{k+1}$ to subproblem \eqref{ADMM1} is expressed explicitly due to its optimality conditions \cite[\S4.2.5]{S.Boyd2011DistributedOptimizationStatistical}, and the $z^{k+1}$ to \eqref{ADMM2} is 
$$\mathop{\text{argmin}}_{z} \Big\{\lambda_{k}\|z\|_{1}+\frac{\rho^{k}}{2}\|u^{k}-z\|^{2}_{2}\Big\},$$
where $u^{k}=-x^{k+1}+\frac{y^{k}}{\rho^{k}}$, its proximal map is the simple and explicit soft threshold operator as follows
$$z^{k+1}_{i}=sign(u^{k}_{i})\Big(|u^{k}_{i}|-\min\{|u^{k}_{i}|,\frac{\lambda_{k}}{\rho^{k}}\}\Big).$$ 

\begin{algorithm}[H]
	\SetAlgoLined
	Initialize $x^{0}$, $z^{0}$, $y^{0}$, $\rho^{0}$, $sn$, $\bar{\epsilon}\in(0,1)$, $\bar{n}\ge 1$, $q$.
	%\KwResult{$x_{k+1}$.}
	
	\For{$k=0,1,2,\ldots$}{
		Compute $x^{k+1}$ by solving \eqref{ADMM1}.\par
		Compute $z^{k+1}$ by solving \eqref{ADMM2}.\par
		Update $y^{k+1}$ by means of \eqref{ADMM3}.\par
		\If{$|i:x^{k+1}_{i}<0|>sn$}{Update $\lambda_{k+1}$ by \eqref{lambda}.}
		% \par
		\eIf{$\mod(k,\bar{n})=1$}{
			$\bar{y}^{k+1}:=y^{k}+\rho^{k}(c-Ax^{k+1}-Bz^{k})$.\par
			Compute spectral step sizes $\alpha_{k}, \beta_{k}$ according to \eqref{equ:alphaRBB} and \eqref{equ:betaRBB} using  parameter \eqref{tauk}.\par
			Compute correlations $\alpha^{corr}, \beta^{corr}$ by \eqref{corr}.\par
			Update penalty parameter $\rho^{k+1}$ according to \eqref{rhoRBB}. 	
		}{
			$\rho^{k+1}=\rho^{k}$.
		}
	}
	\caption{A general scheme for ADMM with adaptive penalty parameter}
	\label{alg:AADMM}
\end{algorithm}

\section{Numerical Experiments}\label{numerical}

%References
\bibliographystyle{plain}
\bibliography{ADMM}

\end{document}